\documentclass[12pt]{amsart}
\usepackage{tipa}
\usepackage{amssymb}
\usepackage{graphicx}
\theoremstyle{definition}

\theoremstyle{remark}

\numberwithin{equation}{section}

\begin{document}
\title[]{{\Large \bf Growth estimates for a class of subharmonic Functions
in a Half Plane  $^{\ast}$ } }
\author{   Pan Guoshuang$^{1,2}$ and Deng Guantie$^{1, \ast\ast}$}%
\address{$^{1}$ Sch. Math. Sci. \& Lab. Math. Com. Sys. \\
   Beijing Normal University\\
    100875  Beijing, The People's Republic of China
    }%
\address{$^{2}$ Department of Public Basic Courses \\
   Beijing  Institute of Fashion and Technology \\
    100029  Beijing, The People's Republic of China
    }%
\email{denggt@bnu.edu.cn}%

\thanks{$\ast $ The project supported by NSFC (Grant No.10671022) and by RFDP (Grant
No.20060027023)
\endgraf $\ast\ast$:Corresponding author.}

\begin{abstract}
 \ \ {  A class of subharmonic  functions represented by
 the modified  kernels are proved to have the  growth estimates
 $u(z)= o(y^{1-\alpha}|z|^{m+\alpha})$ at infinity in the upper half plane  ${\bf
C}_{+}$, which generalizes the growth properties of analytic
functions and harmonic functions. }
 \end{abstract}
\maketitle


 \section*{ 1. Introduction and Main Theorem}

\vspace{0.3cm}

  { Let ${\bf C} $  denote the complex plane
with points $z=x+iy$, where $x, y \in {\bf R}$.  The boundary and
closure of an open $\Omega$ of ${\bf C}$ are denoted by
$\partial{\Omega}$
 and $\overline{\Omega}$ respectively.
 The upper half-plane ${\bf C}_{+}$ is the set
 ${\bf C}_{+}=\{z=x+iy\in {\bf C}:\; y>0\}$, whose boundary is
 $\partial{{\bf C}_{+}}$ .
    We  write $B(z,\rho)$ and $\partial B(z,\rho) $ for the open ball
    and the sphere of radius $\rho$  centered at $z$ in ${\bf C}$.
 We identify $\partial {{\bf C}_{+}}$ with ${\bf R}$.

  For $z\in{\bf C}\backslash\{0\}$, let([\textbf{3}])
$$
 E(z)=(2\pi)^{-1}\log|z|
$$
where $|z|$ is the Euclidean norm . We know that $E$ is
locally integrable in ${\bf C} $.\\
 We define the Green function $G(z,\zeta)$ for the upper half plane
 ${\bf C}_{+}$ by([\textbf{3}])
$$
 G(z,\zeta)=E(z-\zeta)-E(z-\overline{\zeta})
 \qquad z,\zeta\in\overline{{\bf C}_{+}} ,\  z\neq \zeta.\eqno{(1.1)}
$$
We define the Poisson kernel $P(z,\xi)$ when $z\in {\bf C}_{+}$ and
$\xi\in
\partial {\bf C}_{+} $ by
$$
 P(z,\xi)=-\frac{\partial G(z,\zeta)}{\partial
 \eta}\bigg|_{\eta=0}=\frac{y}{\pi|z-\xi|^2}.
$$

  The  Dirichlet problem of upper half plane is to find a function
 $u$ satisfying
$$
 u\in C^2({\bf C}_{+}),
$$
$$
 \Delta u=0,   z\in {\bf C}_{+},
$$
$$
 lim_{z\rightarrow x}u(z)=f(x)\ {\rm nontangentially  \  a.e.}x\in \partial {\bf C}_{+},
$$
where $f$ is a measurable function of ${\bf R}$. The Poisson
integral of the upper half plane is defined by
$$
v(z)=P[f](z)=\int_{{\bf R}}P(z,\xi)f(\xi)d\xi. \eqno{(1.2)}
$$
 We have know that, the Poisson integral $P[f]$ exists if
$$
\int_{{\bf R}}\frac{|f(\xi)|}{1+|\xi|^2} d\xi<\infty. \eqno{(1.3)}
$$
(see [\textbf{4}] and [\textbf{5}])In this paper, we will consider
measurable functions $f$ in ${\bf R}$  satisfying
$$
\int_{{\bf R}}\frac{|f(\xi)|}{1+|\xi|^{2+m}}
d\xi<\infty,\eqno{(1.4)}
$$
where $m$ is a natural number. To obtain a solution of Dirichlet
problem for the boundary date $f$, we use the following modified
functions defined by
$$
 E_m(z-\zeta)=\left\{\begin{array}{ll}
 E(z-\zeta)  &   \mbox{when }   |\zeta|\leq 1,  \\
 E(z-\zeta) -
 \frac{1}{2\pi}\Re(\log\zeta-\sum_{k=1}^{m-1}\frac{z^k}{k\zeta^{k}})
  & \mbox{when}\   |\zeta|> 1.
 \end{array}\right.
$$
Then we can define modified Green function $G_m(z,\zeta)$ and the
modified Poisson
 kernel $P_m(z,\xi)$ by
$$
 G_m(z,\zeta)=E_{m+1}(z-\zeta)-E_{m+1}(z-\overline{\zeta}) \qquad z,\zeta\in\overline{{\bf C}_{+}}, \ z\neq
 \zeta;\eqno{(1.5)}
$$
$$
 P_m(z,\xi)=\left\{\begin{array}{ll}
 P(z,\xi)  &   \mbox{when }   |\xi|\leq 1  ,\\
 P(z,\xi) - \frac{1}{\pi}\Im\sum_{k=0}^{m}\frac{z^k}{\xi^{1+k}}&
\mbox{when}\   |\xi|> 1.
 \end{array}\right.\eqno{(1.6)}
$$
where
 $ z=x+iy, \zeta=\xi+i\eta$.

   Hayman([\textbf{1}]) has proved the
  following result:

\vspace{0.2cm}
 \noindent
{\bf Theorem A } Let $f$ be a measurable function in ${\bf R}$
satisfying (1.3), let $\mu$ be a Borel positive  measure satisfying
$$
\int_{{\bf C}_{+}}\frac{\eta}{1+|\zeta|^{2}} d\mu(\zeta)<\infty.
$$
Write the subharmonic function
$$
u(z)= v(z)+h(z), \quad z\in {\bf C}_{+}
$$
where $v(z)$ be the harmonic function defined by (1.2), $h(z)$ is
defined by
$$
h(z)= \int_{{\bf C}_{+}} G(z,\zeta)d\mu(\zeta)
$$
and $G(z,\zeta)$ is defined by (1.1). Then there exists $z_j\in {\bf
C}_{+},\ \rho_j>0,$ such that
$$
\sum _{j=1}^{\infty}\frac{\rho_j}{|z_j|}<\infty
$$
holds and
$$
u(z)= o(|z|)  \quad  {\rm as}  \ |z|\rightarrow\infty
 $$
holds in ${\bf C}_{+}-G$. where $ G=\bigcup_{j=1}^\infty
B(z_j,\rho_j)$.

  Our aim in this paper is to establish the following theorems.

\vspace{0.2cm}
 \noindent
{\bf Theorem 1} Let $f$ be a measurable function in ${\bf R}$
satisfying (1.4), and $0< \alpha\leq 2$. Let $v(z)$ be the harmonic
function defined by
$$
v(z)= \int_{{\bf R}}P_m(z,\xi)f(\xi)d\xi \quad z\in {\bf C}_{+}
 \eqno{(1.7)}
$$
where $P_m(z,\xi)$ is defined by (1.6). Then there exists $z_j\in
{\bf C}_{+},\ \rho_j>0,$ such that
$$
\sum
_{j=1}^{\infty}\frac{\rho_j^{2-\alpha}}{|z_j|^{2-\alpha}}<\infty
\eqno{(1.8)}
$$
holds and
$$
v(z)= o(y^{1-\alpha}|z|^{m+\alpha})  \quad  {\rm as}  \
|z|\rightarrow\infty   \eqno{(1.9)}
 $$
holds in ${\bf C}_{+}-G$. where $ G=\bigcup_{j=1}^\infty
B(z_j,\rho_j)$.

\vspace{0.2cm}
 \noindent
 {\bf Remark 1} If $\alpha=2$, then (1.8) is a finite sum, the set $G$ is a bounded set,
  so (1.9) holds in ${\bf C}_{+}$.

  Next, we will generalize Theorem 1 to subharmonic functions.

\vspace{0.2cm}
 \noindent
{\bf Theorem 2 } Let $f$ be a measurable function in ${\bf R}$
satisfying (1.4), let $\mu$ be a Borel positive  measure satisfying
$$
\int_{{\bf C}_{+}}\frac{\eta}{1+|\zeta|^{2+m}} d\mu(\zeta)<\infty.
$$
Write the subharmonic function
$$
u(z)= v(z)+h(z), \quad z\in {\bf C}_{+}
$$
where $v(z)$ be the harmonic function defined by (1.7), $h(z)$ is
defined by
$$
h(z)= \int_{{\bf C}_{+}} G_m(z,\zeta)d\mu(\zeta)
$$
and $G_m(z,\zeta)$ is defined by (1.5). Then there exists $z_j\in
{\bf C}_{+},\ \rho_j>0,$ such that (1.8) holds and
$$
u(z)= o(y^{1-\alpha}|z|^{m+\alpha})  \quad  {\rm as}  \
|z|\rightarrow\infty  \eqno{(1.10)}
$$
holds in ${\bf C}_{+}-G$. where $ G=\bigcup_{j=1}^\infty
B(z_j,\rho_j)$ and $0< \alpha<2$.

\vspace{0.2cm}
 \noindent
 {\bf Remark 2} If $\alpha=1, m=0$, this is just the result of
Hamman, so our result (1.10) is the generalization of Theorem A.

\vspace{0.4cm}

\section*{2.   Proof of Theorem }

\vspace{0.3cm}

Let $\mu$ be a positive Borel measure  in ${\bf C},\ \beta\geq0$,
the maximal function $M(d\mu)(z)$ of order $\beta$ is defined by
$$
M(d\mu)(z)=\sup_{ 0<r<\infty}\frac{\mu(B(z,r))}{r^\beta},
$$
then the maximal function $M(d\mu)(z):{\bf C} \rightarrow
[0,\infty)$ is semicontinuous, hence measurable. To see this,
$\forall \lambda >0 $, let $D(\lambda)=\{z\in{\bf
C}:M(d\mu)(z)>\lambda\}$. Fix $z \in D(\lambda)$, then $\exists$
 $r>0$ such that $\mu(B(z,r))>tr^\beta$ for some $t>\lambda$, and
$\exists \delta>0$ satisfying
$(r+\delta)^\beta<\frac{tr^\beta}{\lambda}$. If $|\zeta-z|<\delta$,
then $B(\zeta,r+\delta)\supset B(z,r)$, therefore
$\mu(B(\zeta,r+\delta))\geq
tr^\beta=t(\frac{r}{r+\delta})^\beta(r+\delta)^\beta>\lambda(r+\delta)^\beta$.
Thus $B(z,\delta)\subset D(\lambda)$. This proves that $D(\lambda)$
is open for each $\lambda>0$.

 In order to obtain the result, we
need these lemmas below:

\vspace{0.2cm}
 \noindent
{\bf Lemma 1 } Let $\mu$ be a positive Borel measure in ${\bf C},\
\beta\geq0,\ \mu({\bf C})<\infty, \ \forall \lambda \geq 5^{\beta}
\mu({\bf C})$, set
$$
E(\lambda)=\{z\in{\bf C}:|z|\geq2,M(d\mu)(z) >
\frac{\lambda}{|z|^{\beta}}\}
$$
then $  \exists z_j\in E(\lambda)\  ,\ \rho_j> 0,\ j=1,2,\cdots$,
such that
$$
E(\lambda) \subset \bigcup_{j=1}^\infty B(z_j,\rho_j) \eqno{(2.1)}
$$
and
$$
\sum _{j=1}^{\infty}\frac{\rho_j^{\beta}}{|z_j|^{\beta}}\leq
\frac{3\mu({\bf C})5^{\beta}}{\lambda} .\eqno{(2.2)}
$$
Proof: Let $E_k(\lambda)=\{z\in E(\lambda):2^k\leq |z|<2^{k+1}\}$,
then $ \forall z \in E_k(\lambda), \exists \  r(z)>0$, such that
$\mu(B(z,r(z))) >\lambda(\frac{r(z)}{|z|})^{\beta} $, therefore
$r(z)\leq 2^{k-1}$.
 Since $E_k(\lambda)$ can be covered
by
 the union of a family of balls $\{B(z,r(z)):z\in E_k(\lambda) \}$,
 by the Vitali Lemma([\textbf{2}]), $\exists \  \Lambda_k\subset E_k(\lambda)$,
$\Lambda_k$ is at most countable, such that $\{B(z,r(z)):z\in
\Lambda_k \}$ are disjoint and
$$
E_k(\lambda) \subset
 \cup_{z\in \Lambda_k} B(z,5r(z)),
$$
so
$$
E(\lambda)=\cup_{k=1}^\infty E_k(\lambda) \subset \cup_{k=1}^\infty
\cup_{z\in \Lambda_k} B(z,5r(z)).\eqno{(2.3)}
$$

  On the other hand, note that $ \cup_{z\in \Lambda_k} B(z,r(z)) \subset \{z:2^{k-1}\leq
|z|<2^{k+2}\} $, so that
$$
 \sum_{z \in \Lambda_k}\frac{(5r(z))^{\beta}}{|z|^{\beta}}
\leq 5^\beta\sum_{z\in\Lambda_k}\frac{\mu(B(z,r(z)))}{\lambda} \leq
\frac{5^\beta}{\lambda} \mu\{z:2^{k-1}\leq |z|<2^{k+2}\}.
$$
Hence we obtain
$$
 \sum _{k=1}^{\infty}\sum
_{z \in \Lambda_k}\frac{(5r(z))^{\beta}}{|z|^{\beta}}
 \leq
 \sum _{k=1}^{\infty}\frac{5^\beta}{\lambda} \mu\{z:2^{k-1}\leq |z|<2^{k+2}\}
 \leq
\frac{3\mu({\bf C})5^{\beta}}{\lambda}.\eqno{(2.4)}
$$
  Rearrange $ \{z:z \in \Lambda_k,k=1,2,\cdots\} $ and $
\{5r(z):z \in \Lambda_k,k=1,2,\cdots\}
 $, we get $\{z_j\}$
and $\{\rho_j\}$
 such that (2.1) and
(2.2) hold.

\vspace{0.2cm}
 \noindent
{\bf Lemma 2 } $(1)\  |\Im\sum_{k=0}^{m}\frac{z^k}{\xi^{1+k}}|
\leq\sum_{k=0}^{m-1}\frac{2^ky|z|^k}{|\xi|^{2+k}} ;$\\
$(2)\  |\Im\sum_{k=0}^{\infty}\frac{z^{k+m+1}}{\xi^{k}}|
\leq 2^{m+1}y|z|^m; $\\
$(3)\  |G_m(z,\zeta)-G(z,\zeta)|\leq
\frac{1}{\pi}\sum_{k=1}^{m}\frac{ky\eta|z|^{k-1}}{|\zeta|^{1+k}};$\\
$(4)\  |G_m(z,\zeta)|\leq
\frac{1}{\pi}\sum_{k=m+1}^{\infty}\frac{ky\eta|z|^{k-1}}{|\zeta|^{1+k}}.$\\

  Now we are ready to prove Theorems.

  Throughout the proof, $A$
denote various positive constants.

 \emph{Proof of Theorem 1}

Define the measure $dm(\xi)$ and the kernel $K(z,\xi)$ by
$$
dm(\xi)=\frac{|f(\xi)|}{1+|\xi|^{2+m}} d\xi ,\ \ K(z,\xi)=
P_m(z,\xi)(1+|\xi|^{2+m}).
$$
  For any $\varepsilon >0$, there exists $R_\varepsilon >2$, such that
$$
\int_{|\xi|\geq
R_\varepsilon}dm(\xi)\leq\frac{\varepsilon}{5^{2-\alpha}}.
$$
For every Lebesgue measurable set $E \subset {\bf R}$ , the measure
$m^{(\varepsilon)}$ defined by $m^{(\varepsilon)}(E)
=m(E\cap\{x\in{\bf R}:|x|\geq R_\varepsilon\}) $ satisfies
$m^{(\varepsilon)}({\bf R})\leq\frac{\varepsilon}{5^{2-\alpha}}$,
write
\begin{eqnarray*}
&v_1(z)& =\int_{|\xi-z| \leq 3|z|} P(z,\xi)(1+|\xi|^{2+m})
dm^{(\varepsilon)}(\xi),\\
&v_2(z)&=\int_{|\xi-z| \leq 3|z|}
(P_m(z,\xi)-P(z,\xi))(1+|\xi|^{2+m})
dm^{(\varepsilon)}(\xi), \\
&v_3(z)&=\int_{|\xi-z| > 3|z|} K(z,\xi)dm^{(\varepsilon)}(\xi), \\
&v_4(z)&=\int_{1<|\xi|<R_\varepsilon}K(z,\xi) dm(\xi), \\
&v_5(z)&=\int_{|\xi|\leq1}K(z,\xi) dm(\xi). \\
\end{eqnarray*}
then
$$
|v(z)| \leq |v_1(z)|+|v_2(z)|+|v_3(z)|+|v_4(z)|+|v_5(z)|.
\eqno{(2.5)}
$$
Let $ E_1(\lambda)=\{z\in{\bf C}:|z|\geq2,\exists
t>0,m^{(\varepsilon)}(B(z,t)\cap{\bf R}
)>\lambda(\frac{t}{|z|})^{2-\alpha}\}$, when $ |z|\geq
2R_\varepsilon,\  z \notin E_1(\lambda)
 $, then
$$
\forall t>0,\ m^{(\varepsilon)}(B(z,t)\cap{\bf R}
)\leq\lambda(\frac{t}{|z|})^{2-\alpha}.
$$
  So we have
\begin{eqnarray*}
|v_1(z)|
&\leq& \int_{y\leq|\xi-z| \leq
3|z|}\frac{y}{\pi|z-\xi|^2}2|\xi|^{2+m} dm^{(\varepsilon)}(\xi) \\
&\leq& \frac{2^{2m+5}}{\pi}y|z|^{2+m}\int_{y\leq|\xi-z| \leq
3|z|}\frac{1}{|z-\xi|^2} dm^{(\varepsilon)}(\xi) \\
&=& \frac{2^{2m+5}}{\pi}y|z|^{m+2}\int_{y}^{3|z|} \frac{1}{t^2}
dm_z^{(\varepsilon)}(t).
\end{eqnarray*}
where  $m_z^{(\varepsilon)}(t)=\int_{|\xi-z| \leq t}
dm^{(\varepsilon)}(\xi)$, since for $z \notin E_1(\lambda)$,
\begin{eqnarray*}
\int_{y}^{3|z|} \frac{1}{t^2} dm_z^{(\varepsilon)}(t)
&\leq&  \frac{m_z^{(\varepsilon)}(3|z|)}{(3|z|)^2}+2
\int_{y}^{3|z|} \frac{m_z^{(\varepsilon)}(t)}{t^{3}} dt \\
&\leq& \frac{\lambda}{3^\alpha |z|^2}+2 \int_{y}^{3|z|}
\frac{\lambda\frac{t^{2-\alpha}}{|z|^{2-\alpha}}}{t^{3}} dt \\
&\leq& \frac{\lambda}{ |z|^2}\bigg[\frac{1}{3^\alpha}+
\frac{2}{\alpha}\frac{|z|^\alpha}{y^\alpha}\bigg],
\end{eqnarray*}
so that
$$
|v_1(z)|\leq \frac{2^{2m+5}}{\pi}
\bigg(\frac{1}{3^\alpha}+\frac{2}{\alpha}\bigg)\lambda
y^{1-\alpha}|z|^{m+\alpha}. \eqno{(2.6)}
$$

   By (1) of Lemma 2, we obtain
\begin{eqnarray*}
|v_2(z)|
&\leq& \int_{y\leq|\xi-z| \leq 3|z|}
\frac{1}{\pi}\sum_{k=0}^{m-1}\frac{2^ky|z|^k}{|\xi|^{2+k}}\cdot2|\xi|^{2+m} dm^{(\varepsilon)}(\xi) \\
&\leq&\int_{y\leq|\xi-z| \leq 3|z|}
\sum_{k=0}^{m-1}\frac{2^{k+1}y|z|^k}{\pi}(4|z|)^{m-k} dm^{(\varepsilon)}(\xi) \\
&\leq& \frac{2^{2m+1}}{\pi}\sum_{k=0}^{m-1}\frac{1}{2^k}\frac{1}{5^{2-\alpha}}\varepsilon y|z|^m \\
&\leq& \frac{4^{m-1+\alpha}}{\pi}\varepsilon y|z|^m. \hspace{70mm}
(2.7)
\end{eqnarray*}

  By (2) of Lemma 2, we see that([\textbf{6}])
\begin{eqnarray*}
|v_3(z)|
&\leq& \int_{|\xi-z| > 3|z|}
\bigg|\Im\sum_{k=m}^{\infty}\frac{z^{k+1}}{\pi\xi^{2+k}}\bigg|\cdot2|\xi|^{2+m} dm^{(\varepsilon)}(\xi) \\
&=& \int_{|\xi-z| > 3|z|}
\frac{2}{\pi}\bigg|\Im\sum_{k=0}^{\infty}\frac{z^{k+m+1}}{\xi^{k}}\bigg| dm^{(\varepsilon)}(\xi) \\
&\leq& \frac{2^{m+2}}{\pi}\frac{\varepsilon}{5^{2-\alpha}}
y|z|^m \\
&\leq& \frac{2^{m-2+2\alpha}}{\pi}\varepsilon y|z|^m . \hspace{69mm}
(2.8)
\end{eqnarray*}

  Write
\begin{eqnarray*}
v_4(z)
&=& \int_{1<|\xi|<R_\varepsilon}[P(z,\xi) -\frac{1}{\pi}\Im
\sum_{k=0}^{m}\frac{z^k}{\xi^{1+k}}](1+|\xi|^{2+m}) dm(\xi) \\
&=& v_{41}(z)-v_{42}(z),
\end{eqnarray*}
then
\begin{eqnarray*}
|v_{41}(z)|
&\leq& \int_{1<|\xi|<R_\varepsilon}\frac{y}{\pi|z-\xi|^2}
2|\xi|^{2+m} dm(\xi) \\
&\leq& \frac{2R_\varepsilon^{2+m}y}{\pi}
\int_{1<|\xi|<R_\varepsilon}\frac{1}{(\frac{|z|}{2})^2}
 dm(\xi) \\
&\leq& \frac{2^3R_\varepsilon^{2+m}m({\bf
R})}{\pi}\frac{y}{|z|^2}.\hspace{62mm} (2.9)
\end{eqnarray*}
by (1) of Lemma 2, we obtain
\begin{eqnarray*}
|v_{42}(z)|
&\leq& \int_{1<|\xi|<R_\varepsilon}
\frac{1}{\pi}\sum_{k=0}^{m-1}\frac{2^ky|z|^k}{|\xi|^{2+k}}\cdot 2|\xi|^{2+m} dm(\xi) \\
&\leq& \sum_{k=0}^{m-1}\frac{2^{k+1}}{\pi}y|z|^kR_\varepsilon^{m-k}
m({\bf
R}) \\
&\leq& \frac{2^{m+1}R_\varepsilon^m m({\bf R})}{\pi} y|z|^{m-1}
.\hspace{53mm} (2.10)
\end{eqnarray*}
  In case $|\xi|\leq 1$, note that
$$
K(z,\xi)=P_m(z,\xi)(1+|\xi|^{2+m}) \leq\frac{2y}{\pi|z-\xi|^2},
$$
so that
$$
|v_5(z)|\leq \int_{|\xi|\leq1}\frac{2y}{\pi(\frac{|z|}{2})^2}
dm(\xi)\leq \frac{2^3m({\bf R})}{\pi}\frac{y}{|z|^2}.\eqno{(2.11)}
$$

  Thus, by collecting (2.5), (2.6), (2.7), (2.8), (2.9), (2.10)and
(2.11), there exists a positive constant $A$ independent of
$\varepsilon$, such that if $ |z|\geq 2R_\varepsilon$ and $\ z
\notin E_1(\varepsilon)$, we have
$$
|v(z)|\leq A\varepsilon y^{1-\alpha}|z|^{m+\alpha}.
$$

 Let $\mu_\varepsilon$ be a measure in ${\bf C}$ defined by
$ \mu_\varepsilon(E)= m^{(\varepsilon)}(E\cap{\bf R})$ for every
measurable set $E$ in ${\bf C}$. Take
$\varepsilon=\varepsilon_p=\frac{1}{2^{p+2}}, p=1,2,3,\cdots$, then
there exists a sequence $ \{R_p\}$: $1=R_0<R_1<R_2<\cdots$ such that
$$
\mu_{\varepsilon_p}({\bf C})=\int_{|\xi|\geq
R_p}dm(\xi)<\frac{\varepsilon_p}{5^{2-\alpha}}.
$$
Take $\lambda=3\cdot5^{2-\alpha}\cdot2^p\mu_{\varepsilon_p}({\bf
C})$ in Lemma 1, then $\exists \ z_{j,p}$ and $ \rho_{j,p}$, where
$R_{p-1}\leq |z_{j,p}|<R_p$ such that
$$
\sum _{j=1}^{\infty}(\frac{\rho_{j,p}}{|z_{j,p}|})^{2-\alpha} \leq
\frac{1}{2^{p}}.
$$
So if $R_{p-1}\leq |z|<R_p,\ z\notin G_p=\cup_{j=1}^\infty
B(z_{j,p},\rho_{j,p})$, we have
$$
|v(z)|\leq A\varepsilon_py^{1-\alpha}|z|^{m+\alpha},
$$
thereby
$$
\sum _{p=1}^{\infty}
\sum_{j=1}^{\infty}(\frac{\rho_{j,p}}{|z_{j,p}|})^{2-\alpha} \leq
\sum _{p=1}^{\infty}\frac{1}{2^{p}}=1<\infty.
$$

  Set $ G=\cup_{p=1}^\infty G_p$, then Theorem 1 holds.

 \emph{Proof of Theorem 2}

  Define the measure $dn(\zeta)$ and the kernel $L(z,\zeta)$ by
$$
dn(\zeta)=\frac{\eta d\mu(\zeta)}{1+|\zeta|^{2+m}},\ \
L(z,\zeta)=G_m(z,\zeta)\frac{1+|\zeta|^{2+m}}{\eta}.
$$
then the function $h(z)$ can be written as
$$
h(z)=\int_{{\bf C}_{+}} L(z,\zeta) dn(\zeta).
$$

  For any $\varepsilon >0$, there exists $R_\varepsilon >2$, such that
$$
\int_{|\zeta|\geq
R_\varepsilon}dn(\zeta)<\frac{\varepsilon}{5^{2-\alpha}}.
$$
For every Lebesgue measurable set $E \subset {\bf C}$ , the measure
$n^{(\varepsilon)}$ defined by $n^{(\varepsilon)}(E)
=n(E\cap\{\zeta\in {\bf C}_{+}:|\zeta|\geq R_\varepsilon\}) $
satisfies $n^{(\varepsilon)}({\bf
C}_{+})\leq\frac{\varepsilon}{5^{2-\alpha}}$, write
\begin{eqnarray*}
&h_1(z)&
=\int_{|\zeta-z|\leq\frac{y}{2}}G(z,\zeta)\frac{1+|\zeta|^{2+m}}{\eta}
dn^{(\varepsilon)}(\zeta), \\
&h_2(z)&=\int_{\frac{y}{2}<|\zeta-z|\leq3|z|}G(z,\zeta)\frac{1+|\zeta|^{2+m}}{\eta}
dn^{(\varepsilon)}(\zeta), \\
&h_3(z)&=\int_{|\zeta-z|\leq3|z|}(G_m(z,\zeta)-G(z,\zeta))\frac{1+|\zeta|^{2+m}}{\eta}
dn^{(\varepsilon)}(\zeta), \\
&h_4(z)&=\int_{|\zeta-z|>3|z|}L(z,\zeta)
dn^{(\varepsilon)}(\zeta), \\
&h_5(z)&=\int_{1<|\zeta|<R_\varepsilon}L(z,\zeta) dn(\zeta), \\
&h_6(z)&=\int_{|\zeta|\leq1}L(z,\zeta) dn(\zeta). \\
\end{eqnarray*}
then
$$
h(z)=h_1(z)+h_2(z)+h_3(z)+h_4(z)+h_5(z)+h_6(z). \eqno{(2.12)}
$$
Let $ E_2(\lambda)=\{z\in{\bf C}:|z|\geq2,\exists
t>0,n^{(\varepsilon)}(B(z,t)\cap {\bf C}_{+}
)>\lambda(\frac{t}{|z|})^{2-\alpha}\}, $ when $ |z|\geq
2R_\varepsilon,\ z\notin E_2(\lambda)
 $, then
$$\forall t>0, \ n^{(\varepsilon)}(B(z,t)\cap {\bf C}_{+}
)\leq\lambda(\frac{t}{|z|})^{2-\alpha}.$$

  So we have
\begin{eqnarray*}
|h_1(z)|
&\leq&
\int_{|\zeta-z|\leq\frac{y}{2}}\frac{1}{2\pi}\log\bigg|\frac{\zeta-\overline{z}}{\zeta-z}\bigg|
\frac{1+|\zeta|^{2+m}}{\eta}dn^{(\varepsilon)}(\zeta) \\
&\leq& \int_{|\zeta-z|\leq\frac{y}{2}}\frac{1}{2\pi}
\log\frac{3y}{|\zeta-z|}\frac{2|\zeta|^{2+m}}{\frac{y}{2}}
dn^{(\varepsilon)}(\zeta) \\
&\leq& \frac{2\times (3/2)^{2+m}
}{\pi}\frac{|z|^{2+m}}{y}\int_{|\zeta-z|\leq\frac{y}{2}}\log\frac{3y}{|\zeta-z|}
dn^{(\varepsilon)}(\zeta) \\
&=& \frac{2\times (3/2)^{2+m}
}{\pi}\frac{|z|^{2+m}}{y}\int_0^\frac{y}{2}
\log\frac{3y}{t} dn_z^{(\varepsilon)}(t)\\
&\leq& \frac{2\times (3/2)^{2+m}
}{\pi}\bigg[\frac{\log6}{2^{2-\alpha}}+
\frac{1}{(2-\alpha)2^{2-\alpha}}\bigg]\lambda
y^{1-\alpha}|z|^{m+\alpha}.\hspace{7mm} (2.13)
\end{eqnarray*}
where $ n_z^{(\varepsilon)}(t)=\int_{|\zeta-z| \leq t}
dn^{(\varepsilon)}(\zeta)$.\\

  Note that
$$
|G(z,\zeta)|=|E(z-\zeta)-E(z-\overline{\zeta})|\leq
\frac{y\eta}{\pi|z-\zeta|^2} \eqno{(2.14)}
$$
then by (2.14), we have
\begin{eqnarray*}
|h_2(z)|
&\leq&
\int_{\frac{y}{2}<|\zeta-z|\leq3|z|}\frac{y\eta}{\pi|z-\zeta|^2}
\frac{2|\zeta|^{2+m}}{\eta}
dn^{(\varepsilon)}(\zeta) \\
&\leq&
\frac{2^{2m+5}}{\pi}y|z|^{2+m}\int_{\frac{y}{2}<|\zeta-z|\leq3|z|}
\frac{1}{|z-\zeta|^2} dn^{(\varepsilon)}(\zeta) \\
&=& \frac{2^{2m+5}}{\pi}y|z|^{2+m}\int_\frac{y}{2}^{3|z|}
\frac{1}{t^2} dn_z^{(\varepsilon)}(t)\\
&\leq& \frac{2^{2m+5}}{\pi}y|z|^{2+m}\frac{\lambda}{
|z|^2}\bigg(\frac{1}{3^\alpha}+
\frac{2^{\alpha+1}}{\alpha}\frac{|z|^\alpha}{y^\alpha}\bigg) \\
&\leq& \frac{2^{2m+5}}{\pi}\bigg(\frac{1}{3^\alpha}+
\frac{2^{\alpha+1}}{\alpha}\bigg)\lambda y^{1-\alpha}|z|^{m+\alpha}.
\hspace{35mm} (2.15)
\end{eqnarray*}

   By (3) of Lemma 2 , we obtain
\begin{eqnarray*}
|h_3(z)|
&\leq&
\int_{|\zeta-z|\leq3|z|}\frac{1}{\pi}\sum_{k=1}^{m}\frac{ky\eta|z|^{k-1}}{|\zeta|^{1+k}}
\frac{2|\zeta|^{2+m}}{\eta} dn^{(\varepsilon)}(\zeta) \\
&\leq&
\int_{|\zeta-z|\leq3|z|}\frac{2}{\pi}\sum_{k=1}^{m}ky|z|^{k-1}(4|z|)^{m-k+1}
 dn^{(\varepsilon)}(\zeta) \\
&\leq&
\frac{2^{2m+1}}{\pi}\sum_{k=1}^{m}\frac{k}{4^{k-1}}\frac{1}{5^{2-\alpha}}\varepsilon y|z|^m \\
&\leq& \frac{2^{2m+2\alpha+1}}{9\pi}\varepsilon y|z|^m.
\hspace{65mm} (2.16)
\end{eqnarray*}

  By (4) of Lemma 2, we see that
\begin{eqnarray*}
|h_4(z)|
&\leq&
\int_{|\zeta-z|>3|z|}\frac{1}{\pi}\sum_{k=m+1}^{\infty}\frac{ky\eta|z|^{k-1}}{|\zeta|^{1+k}}|
\frac{2|\zeta|^{2+m}}{\eta} dn^{(\varepsilon)}(\zeta) \\
&\leq& \int_{|\zeta-z|>3|z|}\frac{2}{\pi}\sum_{k=m+1}^{\infty}ky
\frac{|z|^{k-1}}{(2|z|)^{k-m-1}}
dn^{(\varepsilon)}(\zeta) \\
&\leq& \frac{2^{m+2}}{\pi}\sum_{k=m+1}^{\infty}\frac{k}{2^{k}}
 \frac{1}{5^{2-\alpha}}\varepsilon y|z|^m \\
&\leq& \frac{4^{\alpha-1}(m+2)}{\pi}\varepsilon y|z|^m .
\hspace{59mm} (2.17)
\end{eqnarray*}

  Write
\begin{eqnarray*}
h_5(z)
&=&
\int_{1<|\zeta|<R_\varepsilon}[G(z,\zeta)+(G_m(z,\zeta)-G(z,\zeta))]\frac{1+|\zeta|^{2+m}}{\eta} dn(\zeta) \\
&=& h_{51}(z)+h_{52}(z),
\end{eqnarray*}
then we obtain by (2.14)
\begin{eqnarray*}
|h_{51}(z)|
&\leq& \int_{1<|\zeta|<R_\varepsilon}
\frac{y\eta}{\pi|z-\zeta|^2}\frac{2|\zeta|^{2+m}}{\eta} dn(\zeta) \\
&\leq&
\frac{2R_\varepsilon^{2+m}}{\pi}y\int_{1<|\zeta|<R_\varepsilon}
\frac{1}{(\frac{|z|}{2})^2} dn(\zeta) \\
&\leq& \frac{2^3R_\varepsilon^{2+m}n({\bf
C}_{+})}{\pi}\frac{y}{|z|^2}.\hspace{58mm} (2.18)
\end{eqnarray*}
by (3) of Lemma 2 , we obtain
\begin{eqnarray*}
|h_{52}(z)|
&\leq&
\int_{1<|\zeta|<R_\varepsilon}\frac{1}{\pi}\sum_{k=1}^{m}\frac{ky\eta|z|^{k-1}}{|\zeta|^{1+k}}
\frac{2|\zeta|^{2+m}}{\eta} dn(\zeta) \\
&\leq& \frac{2}{\pi}\sum_{k=1}^{m}ky|z|^{k-1}R_\varepsilon^{m-k+1}
 n({\bf C}_{+}) \\
&\leq& \frac{m(m+1)R_\varepsilon^{m}n({\bf C}_{+}) }{\pi}
y|z|^{m-1}.\hspace{42mm} (2.19)
\end{eqnarray*}

  In case $|\zeta|\leq 1$, by (2.14), we have
$$
|L(z,\zeta)|\leq \frac{y\eta}{\pi|z-\zeta|^2}\frac{2}{\eta}
=\frac{2y}{\pi|z-\zeta|^2},
$$
so that
$$
|h_6(z)|\leq \int_{|\zeta|\leq1}\frac{2y}{\pi(\frac{|z|}{2})^2}
dn(\zeta) \leq \frac{2^3n({\bf C}_{+})}{\pi}\frac{y}{|z|^2}.
\eqno{(2.20)}
$$

  Thus, by collecting (2.12), (2.13), (2.15), (2.16),
(2.17), (2.18), (2.19) and (2.20), there exists a positive constant
$A$ independent of $\varepsilon$, such that if $ |z|\geq
2R_\varepsilon$ and $\  z \notin E_2(\varepsilon)$, we have
$$
 |h(z)|\leq A\varepsilon y^{1-\alpha}|z|^{m+\alpha}.
$$

  Similarly, if $z\notin G$, we have
$$
h(z)= o(y^{1-\alpha}|z|^{m+\alpha})\quad  {\rm as} \
|z|\rightarrow\infty. \eqno{(2.21)}
$$

by (1.9) and  (2.21), we obtain
$$
u(z)=v(z)+h(z)= o(y^{1-\alpha}|z|^{m+\alpha})\quad  {\rm as} \
|z|\rightarrow\infty
$$
hold in  ${\bf C}_{+}-G$.

\begin{center}

\end{center}

\end{document}